\RequirePackage{ifpdf}
\ifpdf 
\documentclass[pdftex]{sigma}
\else
\documentclass{sigma}
\fi

\def\tr{\mathop{\rm tr}}
\def\del{\nabla}
\def\hook{\vrule height 0pt depth 0.4pt width 4pt \vrule height 8pt depth 0.4pt
\kern 3pt}

\def\fpd#1#2{\frac{\partial #1}{\partial #2}}
\def\R{{\rm I\kern-.20em R}}
\def\cinfty#1{C^{\scriptscriptstyle\infty}(#1)}
\def\vectorfields#1{{\cal X}(#1)}

\def\ov#1{\overline{#1}}
\def\T{{\mbox{\bf T}}}
\def\H#1{{#1}^{\scriptscriptstyle H}}
\def\V#1{{#1}^{\scriptscriptstyle V}}
\def\DV#1{\V{\rm D}_{#1}}

\def\dh#1{d^{\scriptscriptstyle H}_{#1}}

\begin{document}


\allowdisplaybreaks

\renewcommand{\PaperNumber}{024}

\FirstPageHeading

\renewcommand{\thefootnote}{$\star$}

\ShortArticleName{First Integrals of the Geodesic Flow of a
Finsler Manifold}

\ArticleName{A Recursive Scheme of First Integrals\\ of the
Geodesic Flow of a Finsler Manifold\footnote{This paper is a
contribution to the Proceedings of the Workshop on Geometric
Aspects of Integ\-rable Systems
 (July 17--19, 2006, University of Coimbra, Portugal).
The full collection is available at
\href{http://www.emis.de/journals/SIGMA/Coimbra2006.html}{http://www.emis.de/journals/SIGMA/Coimbra2006.html}}}

\Author{Willy SARLET}

\AuthorNameForHeading{W.~Sarlet}

\Address{Department of Mathematical Physics and Astronomy,
Ghent University,\\
Krijgslaan 281, B-9000 Ghent, Belgium}
\Email{\href{mailto:willy.sarlet@ugent.be}{willy.sarlet@ugent.be}}
\URLaddress{\url{http://maphyast.ugent.be}}

\ArticleDates{Received October 30, 2006, in f\/inal form January
17, 2007; Published online February 13, 2007}

\Abstract{We review properties of so-called special conformal
Killing tensors on a Riemannian manifold $(Q,g)$ and the way they
give rise to a Poisson--Nijenhuis structure on the tangent bundle
$TQ$. We then address the question of generalizing this concept to
a~Finsler space, where the metric tensor f\/ield comes from a
regular Lagrangian function~$E$, homogeneous of degree two in the
f\/ibre coordinates on $TQ$. It is shown that when a symmetric
type~(1,1) tensor f\/ield $K$ along the tangent bundle projection
$\tau\colon TQ\rightarrow Q$ satisf\/ies a~dif\/ferential
condition which is similar to the def\/ining relation of special
conformal Killing tensors, there exists a direct recursive scheme
again for f\/irst integrals of the geodesic spray. Involutivity of
such integrals, unfortunately, remains an open problem.}


\Classification{37J35; 53C60; 70H06}

\begin{flushright}
\it This paper is based on joint work with Fien Vermeire and Mike
Crampin.
\end{flushright}

\section{Introduction: special conformal Killing tensors}

The work presented here is inspired by the theory and applications
of so-called {\em special conformal Killing tensors\/} (or {\em
Benenti tensors\/}) on a (pseudo-) Riemannian manifold. The study
of possible generalizations to Finsler manifolds should be seen as
a f\/irst step towards further generalizations to arbitrary metric
tensor f\/ields along the tangent bundle projection (coming, for
example, from the Hessian of a regular Lagrangian); the aim in the
end would be to arrive at a scheme which can cope with integrable
or Hamilton--Jacobi separable systems with non-quadratic integrals
in involution, as opposed to the subclass of St\"ackel systems
which is governed by special conformal Killing tensors. We start
by reviewing a number of essential features related to such
tensors.

Let $g$ be a (pseudo-) Riemannian metric on a manifold $Q$.

\begin{definition} A symmetric tensor $J_{ij}$ determines a {\em
special conformal Killing tensor\/} (or Benenti tensor) w.r.t.\
$g$ if its covariant derivatives with respect to the Levi-Civita
connection satisfy,
\begin{gather*}
J_{ij|k} = \frac{1}{2} (\alpha_i g_{jk} + \alpha_j g_{ik}),
\end{gather*}
for some 1-form $\alpha$.
\end{definition}
It follows that
\begin{gather*}
\sum_{ijk} J_{ij|k} = \sum_{ijk} \alpha_i g_{jk},
\end{gather*}
i.e.\ $J$ is a conformal Killing tensor. Further properties
(referring to $J$ as (1,1)-tensor: $J^i_j=g^{ik}J_{kj}$) are that
$J$ is a conformal Killing tensor of gradient type,
\begin{gather*}
\alpha=df=d(\tr J),
\end{gather*}
and has vanishing Nijenhuis torsion $N_J$, i.e.
\begin{gather*}
{}[JX,JY]+J^2([X,Y]) - J([JX,Y]+[X,JY])=0 \qquad \forall\;
X,Y\in\vectorfields{Q}.
\end{gather*}
Moreover, assuming $J$ is non-singular, its cofactor tensor $A$:
\begin{gather*}
A\,J=(\det J) I
\end{gather*}
is itself a Killing tensor, so that $F=(1/2) A^{ij}(q)p_ip_j$ is a
quadratic f\/irst integral of the system with Hamiltonian $H=(1/2)
g^{ij}(q)p_ip_j$ on $T^*Q$.

As said above, metrics admitting a special conformal Killing
tensor determine a subclass of St\"ackel systems, i.e.\ orthogonal
separable systems in the sense of Hamilton--Jacobi with $n$
quadratic integrals in involution. We brief\/ly sketch how this
works, beginning with the integrability structures on $T^*Q$ and
$T^*Q\times\mathbb{R}$, determined by $J$.

On $T^*Q$, consider $\widetilde J$, the {\it complete lift\/} of
$J$, which in coordinates is given by
\begin{gather*}
\widetilde{J}=J^i_j\left(\fpd{}{q^i}\otimes dq^j +
\fpd{}{p_j}\otimes dp_i\right) +
p_k\left(\fpd{J^k_i}{q^j}-\fpd{J^k_j}{q^i}\right)
\fpd{}{p_i}\otimes dq^j.
\end{gather*}
If $\omega=d\theta$ denotes the standard symplectic structure and
\begin{gather*}
\omega_1= d(J\theta) = d(J^i_jp_idq^j),
\end{gather*}
we have that
\begin{gather*}
i_{\widetilde{J}(X)}\omega = i_X\omega_1 \qquad \forall \;
X\in\vectorfields{T^*Q}.
\end{gather*}
It follows that $\widetilde{J}$ is symmetric with respect to
$\omega$. Moreover, we have $d\omega_1=0$ obviously, and
$N_{\widetilde{J}}=0$ as a result of $N_J=0$. These are the
properties which are suf\/f\/icient to ensure that a~special
conformal Killing tensor $J$ on $(Q,g)$ determines a
Poisson--Nijenhuis structure on $T^*Q$ with $\widetilde{J}$ as
recursion operator. In addition, we have the interesting relations
\begin{gather*}
dd_{\widetilde J}h=d(\tr J)\wedge dh, \qquad dd_{\widetilde J}
(\tr J)=0,
\end{gather*}
which give rise to an extension of the two compatible Poisson
structures to $T^*Q\times\mathbb{R}$. Details about all of this
can be found in \cite{CST}.

Observe that the special conformal Killing tensor $J$ plays a kind
of double role for the complete integrability of $H=(1/2)
g^{ij}p_ip_j$. First of all, we have that $J+sI$ is a special
conformal Killing tensor for all $s$, and thus gives rise to a
family of corresponding Killing tensors $A(s)$, and in this way to
a hierarchy of quadratic integrals. Secondly, the double Poisson
structure created by $J$ truly helps to show that such integrals
are in involution, in fact with respect to both brackets (although
other, more direct techniques for showing involutivity may
sometimes be successful as well, see e.g.\ \cite{TopMat}).
Incidentally, the reference just cited also addresses another
important aspect of the integrals, namely their functional
independence, for which it is shown to be enough that the
eigenvalues of $J$ are dif\/ferent at one point. To complete the
construction of St\"ackel systems then, one can proceed as
follows. For a function $V$ on $Q$ to be an admissible potential,
it is necessary and suf\/f\/icient that $V$ satisf\/ies
\begin{gather*}
dd_{\widetilde J}V = d(\tr J)\wedge dV.
\end{gather*}
Admissibility of $V$ means that
\begin{gather*}
h=\frac{1}{2} g^{ij}p_ip_j + V
\end{gather*}
is still separable, and there are then suitable modif\/ications of
the f\/irst integrals of the kinetic energy part which constitute
an amended set of quadratic integrals in involution. To avoid
forgetting a number of interesting references concerning details
and background about these results, we limit ourselves to citing
the excellent review paper \cite{Benenti} and refer to the list of
references therein.

\section{A tangent bundle view of special conformal Killing tensors}

We explain some of the results from \cite{SarVer} here;
understanding a Lagrangian perspective of what precedes is
essential for the generalization to Finsler spaces we have in
mind.

Let $(Q,g)$ as before be a pseudo-Riemannian manifold and consider
the function
\begin{gather*}
L= \frac{1}{2} g_{ij}(q)u^iu^j \in\cinfty{TQ}.
\end{gather*}
Then, $TQ$ becomes a symplectic manifold via the
Poincar\'e--Cartan 2-form
\begin{gather*}
\omega_L=d\theta_L, \qquad \mbox{with}\qquad
\theta_L=S(dL)=\fpd{L}{u^i}dq^i,
\end{gather*}
where $S$ is the canonically def\/ined `vertical endomorphism' on
$TQ$.

Consider then again a type (1,1) tensor f\/ield $J$ on $Q$. The
following is an interesting result of \cite{SarVer}: it gives a
concise, intrinsic characterization of special conformal Killing
tensors, directly in their type (1,1) appearance.

\begin{theorem} $J$ is a special conformal Killing tensor if and
only if it is symmetric with respect to $g$ and  satisfies the
following relation for some function $f\in\cinfty{Q}$:
\begin{gather}
\del J =\frac{1}{2}(\T\otimes \dh{}\!f - X_f\otimes \theta_L).
\label{delJ}
\end{gather}
\end{theorem}

Obviously, we have a bit of work to do here to explain the
symbols in this formulation. In~principle, this whole formula is
about tensor f\/ields along the tangent bundle projection
$\tau\colon TQ\rightarrow Q$, though some of its ingredients so
far happen to live on the base manifold $Q$. The operations which
are being used rely on the (non-linear) connection associated to a
given second-order f\/ield $\Gamma$ (here coming from the
Euler--Lagrange equations of $L$). If we express the coordinate
representation of a general second-order f\/ield $\Gamma$ as
\begin{gather*}
\Gamma = u^i\fpd{}{q^i} + F^i(q,u)\fpd{}{u^i},
\end{gather*}
then the basis of horizontal vector f\/ields $H_i$ on $TQ$,
def\/ining the connection, is given by
\begin{gather*}
H_i = \fpd{}{q^i} - \Gamma^j_i(q,u)\,\fpd{}{u^j},
\qquad\mbox{with}\qquad \Gamma^j_i= - \frac{1}{2}\fpd{F^j}{u^i}.
\end{gather*}
We have
\begin{gather*}
\dh{}f=H_i(f)dq^i \qquad\mbox{and}\qquad X_f\hook g = - \dh{}\!f.
\end{gather*}
$\T$ is the canonical vector f\/ield along $\tau$ (essentially the
identity map on $TQ$), given by
\begin{gather*}
\T = u^i\fpd{}{q^i},
\end{gather*}
and $\del$, f\/inally, is the {\it dynamical covariant
derivative\/} operator. It acts like $\Gamma$ on functions on
$TQ$, and is further determined by
$\del(\partial/\partial{q^i})=\Gamma^j_i\partial/\partial{q^j}$
and by duality on $dq^j$. Up to now, however, $f$ and $J$ live on
$Q$, so that in fact $\dh{}f=df$ and, taking into account that we
are in a situation where $\Gamma^j_i=\Gamma^j_{ik}(q)u^k$, we have
$(\del J)^i_j = J^i_{j|k}u^k$. Observe that it easily follows from
\eqref{delJ}, by taking a trace, that $ f=\tr J$.

If $J$ is an arbitrary (1,1) tensor on $Q$, its vertical lift
$\V{J}$ def\/ines a kind of alternative almost tangent structure
on $TQ$, at least we have ${\V{J}}^2=0$ and $N_{\V{J}}=0$.
Therefore, it is a natural construction to let $\V{J}$ take over
the role of $S$ for the purpose of def\/ining a second 2-form and
the subsequent construction of a type (1,1) tensor $R$ on $TQ$. In
other words, natural tangent bundle constructions lead us to
def\/ine $R$ on $TQ$ by
\begin{gather}
i_{R(X)}d(S(dL))=i_X d(\V{J}(dL)) \qquad \forall \; X\in
\vectorfields{TQ}, \label{R}
\end{gather}
and if $N_J=0$ on $Q$ (as in the case of a special conformal
Killing tensor), we have $N_R=0$ on~$TQ$ as well, so that $R$
becomes the recursion operator of a Poisson--Nijenhuis structure
on~$TQ$. In fact, $R$ is precisely the pullback of $\widetilde J$
under the Legendre transform associated to the given regular
Lagrangian $L$.

\section{Generalization: Finsler manifolds}

Let $(Q,E)$ now be a Finsler space, that is to say, $E\colon
TQ\rightarrow \mathbb{R}$ is the square of a Finsler function~$F$,
is homogeneous of degree two in the $u^i$ ($E(0)=0$), and
$\omega_E=d\theta_E$ is non-degenerate on the slit tangent bundle
$TQ\setminus\{0\}$.

Putting
\begin{gather*}
g_{ij}(q,u)= \fpd{^2E}{u^i\partial u^j},
\end{gather*}
we have a (0,2) tensor f\/ield along the tangent bundle projection
$\tau$, and consider its canonical spray $\Gamma$ on
$TQ\setminus\{0\}$ (i.e.\ the Euler--Lagrange equations of $E$).
Note that in view of the homogeneity,
\begin{gather*}
E=\frac{1}{2} g(\T,\T)= \frac{1}{2} g_{ij}u^iu^j,
\end{gather*}
and we also have
\begin{gather*}
\del g=0, \qquad \del\T=0.
\end{gather*}
It is then easy to deduce the following additional properties:
\begin{gather*}
\del E=\Gamma(E)=0, \qquad \dh{}E=0, \qquad \theta_E=\T\hook g,
\qquad \nabla\theta_E=0.
\end{gather*}

It may be worthwhile to make the following preliminary remark to
avoid confusion. There is often a lively debate  in Finsler
geometry about the selection of the best possible linear
connection to replace the Levi-Civita connection of Riemannian
geometry. Most of the time, about four dif\/ferent possibilities
are considered (see e.g.\ \cite{BCS}), depending on whether one
wants the connection to be as metrical as possible, or whether one
wants to eliminate as much torsion as possible. We do not need to
enter into this debate here, because all we need for deriving our
main results in this section is the dynamical covariant derivative
operator $\nabla$, as it was specif\/ied in the preceding section:
it comes from the non-linear connection associated to any
second-order dif\/ferential equation f\/ield and acts on the full
algebra of tensor f\/ields along the tangent bundle projection
$\tau$.

Inspired by the theorem of the preceding section, assume now that
we have a type (1,1) tensor~$K$ along $\tau$, which is symmetric
with respect to $g$, i.e.\ $g(KX,Y)=g(X,KY)$, and satisf\/ies
a~relation of the form
\begin{gather}
\del K = \frac{1}{2} (\T\otimes \alpha - X_\alpha\otimes
\theta_E), \label{K}
\end{gather}
for some 1-form $\alpha$ along $\tau$, and with $X_\alpha$
def\/ined by $X_\alpha\hook g= -\alpha$. It would be quite natural
within this Finsler environment to assume that $K$ is homogeneous
of degree zero as well, but as we will show now, such an extra
assumption is not even required to arrive at a quite remarkable
hierarchy of f\/irst integrals for the canonical spray $\Gamma$.

An immediate consequence of \eqref{K}, which follows from taking a
trace, is that
\begin{gather*}
\alpha(\T)=\del \tr K.
\end{gather*}
As a result, if we put
\begin{gather*}
h_0 = k_0 =E, \qquad k_1=\frac{1}{2} g(\T,K\T),
\end{gather*}
it is easy to show that also
\begin{gather*}
h_1 = k_1 - E\tr K
\end{gather*}
is a f\/irst integral. This is the start for the hierarchy of
f\/irst integrals, for which we found the following explicit
recursive scheme. Def\/ine
\begin{gather*}
k_j=\frac{1}{2} g(\T,K^j\T), \qquad\mbox{and}\qquad a_j
=\frac{1}{j}\tr K^j,\qquad j=1,2,\ldots .
\end{gather*}
\begin{lemma} From the fundamental assumption \eqref{K} on $K$, it
follows that
\begin{gather}
\del a_j = \alpha(K^j\T), \label{delan}
\\
\del k_j = \sum_{i=1}^j (\del a_i) k_{j-i}. \label{delkn}
\end{gather}
\end{lemma}
\begin{proof} It is easy to show by induction that \eqref{K} implies
that
\begin{gather*}
\del K^{j+1} = \frac{1}{2}\, \sum_{i=0}^j (K^i\alpha\otimes
K^{j-i}\T - K^i\theta_E \otimes K^{j-i}X_\alpha).
\end{gather*}
Taking a trace, \eqref{delan} immediately follows. Moreover, since
\begin{gather*}
\del k_j = \frac{1}{2} g(\T,\del K^j(\T)),
\end{gather*}
a direct computation, using the formula for $\del K^j$ just
derived, leads to \eqref{delkn}. \end{proof}

Next, we introduce auxiliary functions $\phi_k$, def\/ined
recursively by
\begin{gather}
\phi_1=0,\qquad \phi_k= \frac{1}{2}\sum_{i=1}^{k-1} a_ia_{k-i} -
\frac{1}{k}\sum_{i=2}^{k-1}(k-i)\phi_i a_{k-i},\qquad k\geq 2.
\label{phi}
\end{gather}

\begin{lemma} We have\ $\del\phi_j= \sum\limits_{l=1}^{j-1}a_l\del a_{j-l}
- \sum\limits_{l=2}^{j-1}\phi_l\del a_{j-l}, \quad j\geq 2$.
\end{lemma}

\begin{proof} From $\phi_2=(1/2) a_1^2$, we get
$\del\phi_2=a_1\del a_1$ and the property is clearly true for
$j=2$. Assume it is valid for all $j$ up to $k-1$ and now act with
$\del$ on $\phi_k$ as def\/ined by \eqref{phi}. We get
\begin{gather*}
\del\phi_k = \sum_{i=1}^{k-1}a_i\del a_{k-i} - \frac{1}{k}
\sum_{i=2}^{k-1}(k-i)\phi_i\del a_{k-i}
-\frac{1}{k}\sum_{i=2}^{k-1}(k-i)a_{k-i}\sum_{l=1}^{i-1}a_l\del
a_{i-l}
\\ \phantom{\del\phi_k =}{}
{} +\frac{1}{k}\sum_{i=2}^{k-1}(k-i)a_{k-i}
\sum_{l=2}^{i-1}\phi_l\del a_{i-l}.
\end{gather*}
In order to collect coef\/f\/icients of $\del a_j$, we put $j=k-i$
in the second term on the right and $j=i-l$ in the last two terms
(and adjust the summations). Moreover, we subsequently interchange
the two summations in those last two terms and split of\/f one
term to have the common upper bound $k-3$ for $j$. The result
reads
\begin{gather*}
\del\phi_k = \sum_{i=1}^{k-1}a_i\del a_{k-i} - \frac{1}{k}
\sum_{j=1}^{k-2}j\,\phi_{k-j}\del a_j
\\ \phantom{\del\phi_k =}
{} -\frac{1}{k}\sum_{j=1}^{k-3}\del
a_j\Biggl(\sum_{i=j+1}^{k-1}(k-i)a_{k-i}a_{i-j} -
\sum_{i=j+2}^{k-1}(k-i)a_{k-i}\phi_{i-j}\Biggr) -
\frac{1}{k}\,a_1^2\del a_{k-2}.
\end{gather*}
In the terms between brackets, we now change $i-j$ to $l$, thus
getting
\begin{gather*}
\Biggl(\sum_{l=1}^{k-j-1}(k-j-l)a_{k-j-l}\,a_l -
\sum_{l=2}^{k-j-1}(k-j-l)a_{k-j-l}\phi_l\Biggr)
\end{gather*}
which upon closer inspection can be seen to be equal to
$(k-j)\phi_{k-j}$. It is now a matter of splitting of\/f the term
for $j=k-2$ in the second sum on the right also, to see a few
cancellations presenting themselves, which then readily produce
the desired result. \end{proof}

Finally, put $b_k=\phi_k-a_k$ and def\/ine
\begin{gather*}
h_l = k_l + \sum_{i=1}^l b_i\,k_{l-i}, \qquad l\geq 1.
\end{gather*}

\begin{theorem} Let $K$ be a type $(1,1)$ tensor along
$\tau\colon TQ\rightarrow Q$, which is symmetric with respect to
the Finsler metric $g$ and satisf\/ies condition \eqref{K} for
some $\alpha$. Then the functions $h_l$ are f\/irst integrals of
the geodesic spray of the Finsler metric $g$, for all $l$.
\end{theorem}

\begin{proof} From Lemma~1 we have, using also $\del k_0=0$,
\begin{gather*}
\del h_l= \sum_{i=1}^l (\del \phi_i)k_{l-i} +
\sum_{i=1}^{l-1}b_i\sum_{j=1}^{l-i}(\del a_j)k_{l-i-j}.
\end{gather*}
Putting $m=i+j$ in the second term and interchanging the double
summation, this becomes
\begin{gather*}
\del h_l= \sum_{i=2}^l \Biggl((\del \phi_i) + \sum_{m=1}^{i-1}b_m
(\del a_{i-m})\Biggr) k_{l-i} + (\del \phi_1)k_{l-1}.
\end{gather*}
This is clearly zero in view of $\phi_1=0$ and Lemma~2.
\end{proof}

It is interesting to verify that one can rewrite the recursive
formula for integrals in the following alternative way:
\begin{gather}
h_l=\frac{1}{2} g(\T,B_{l}\T) \qquad\mbox{with}\qquad B_{l}=b_l\,I
+ B_{l-1}K,\qquad B_{0}=I. \label{hnbis}
\end{gather}
Moreover, the $b_l$ turn out to be the coef\/f\/icients of the
characteristic polynomial of $K$: indeed, if $n$ is the dimension
of $Q$, one can show that
\begin{gather}
\det(\lambda\,I-K)= \sum_{i=0}^nb_i\lambda^{n-i} \qquad (b_0=1),
\label{char}
\end{gather}
which means that in particular $b_n=(-1)^n\det K$.

Obviously, this whole construction applies in particular to the
Riemannian case with the assumption that $K$ is a special
conformal Killing tensor on $Q$; in fact, the expression
\eqref{hnbis} looks like a direct generalization of a result of
Benenti (see the f\/irst theorem in section~7 of \cite{Benenti}).

The parallel with the Riemannian case goes further, in the sense
that, if $K$ is assumed to be non-singular, there is again a
f\/irst integral associated to the cofactor tensor. Observe, by
the way, that the quadratic integral $F=(1/2) A^{ij}(q)p_ip_j$
referred to in section~1 for the Riemannian case, can be written
in a tangent bundle set-up in the coordinate free format $F=(1/2)
g(A\T,\T)$, and this statement happens to extend to the Finslerian
situation.

\begin{theorem} Let $g$ be a Finsler metric, and $K$ a symmetric
$(1,1)$ tensor satisfying condition~\eqref{K}. Then, if $K$ is
non-singular and $A$ is its cofactor tensor, the function
$g(A\T,\T)$ is a f\/irst integral of the geodesic spray.
\end{theorem}

\begin{proof} We have that $\del A=\del(\det K)\!K^{-1} - (\det K)
K^{-1}\del K\,K^{-1}$ and it is a general property for  derivation
of degree zero such as $\del$ that $\del\log(\det K)=
\tr(K^{-1}\del K)$. It follows that
\begin{gather*}
\del(g(A\T,\T)) = g(\del A(\T)T)
\\ \phantom{\del(g(A\T,\T))}
 {}= (\det K)(\tr(K^{-1}\del K)g(K^{-1}\T,\T) - g(\del KK^{-1}(\T),K^{-1}(\T))).
\end{gather*}
From the assumption \eqref{K} on $K$, it readily follows that
$\tr(K^{-1}\del K)=\alpha(K^{-1}(\T))$ and
\begin{gather*}
g(\del K\,K^{-1}(\T), K^{-1}(\T)) =
\alpha(K^{-1}(\T))g(\T,K^{-1}\T).
\end{gather*}
Hence, $\del(g(A\T,\T))=0$.
\end{proof}

The hierarchy of f\/irst integrals derived above can now be
established via the cofactor technique as well. Indeed, if $K$
satisf\/ies the fundamental condition \eqref{K}, then so does
$K+sI$ for all real values $s$. Its cofactor
$A(s)=\sum\limits_{j=0}^{n-1}A_{j+1} s^{j}$ therefore also gives
rise to a hierarchy of f\/irst integrals.

In fact, we have
\begin{gather*}
(K+sI)A(s)=\sum_{j=0}^{n-1}KA_{j+1}s^j + \sum_{j=1}^nA_js^j.
\end{gather*}
On the other hand, using \eqref{char} we get
\begin{gather*}
\det(K+sI)= (-1)^n\det((-s)I-K)=\sum_{i=0}^n(-1)^ib_is^{n-i}.
\end{gather*}
It follows that
\begin{gather*}
A_n=I, \qquad \mbox{and}\qquad KA_1=(-1)^nb_nI,
\end{gather*}
which shows that $A_1=A$, while the other coef\/f\/icients in the
expansion of $A(s)$ have to satisfy the recursive relation
\begin{gather*}
KA_{j+1} + A_j = (-1)^{n-j}b_{n-j} I, \qquad j=1,\ldots,n-1.
\end{gather*}
Going down this scheme from top to bottom, it is easy to see that
$A_{n-j}= (-1)^jB_j$, which shows (see \eqref{hnbis}) that the
coef\/f\/icients in the expansion of $A(s)$ determine the same
f\/irst integrals $h_l$ up to sign. Notice that, in particular,
$A=A_1=(-1)^{n-1}B_{n-1}$, from which it follows that
\begin{gather*}
B_n = B_{n-1}K + b_nI= (-1)^{n-1}AK + (-1)^n(\det K)I =0.
\end{gather*}
One can verify also that $b_{n+1}\equiv 0$. Hence, the sequence of
f\/irst integrals actually terminates with $h_n=0$ and the
cofactor tensor $A$ of $K$ determines the last non-zero integral
$h_{n-1}$ in the hierarchy.

\section{The problem of involutivity}

We have not succeeded in coming anywhere near to proving that the
integrals we obtained in the previous section would be in
involution with respect to the Poisson bracket associated to the
symplectic form $d\omega_E$. One would expect to have better
chances if one could use the tensor~$K$ for constructing a
recursion tensor $R$ of a Poisson--Nijenhuis structure on $TQ$,
just as in the Riemannian case. Observe that the construction
\eqref{R} of $R$ in the Riemannian case in fact makes sense under
much more general circumstances. Indeed, the vertical lift $\V{J}$
remains well def\/ined if we replace the basic tensor $J$ on $Q$
by a tensor $J$ along the projection $\tau$, and the function $L$
can be taken to be any regular Lagrangian on a tangent bundle
$TQ$. Such an $R$ then is determined, through its action on
vertical and horizontal lifts of vector f\/ields $X$ along $\tau$,
by formulas of the form (see \cite{VSC})
\begin{gather*}
R(\V{X}) = \V{(\ov{K}X)},
\\[1ex]
R(\H{X}) = \H{(KX)} + \V{(UX)},
\end{gather*}
where $\ov{K}$ is the transpose of $K$ (with respect to $g$) and
$K$ and $U$ are def\/ined by
\begin{gather*}
g(KX,Y) = \DV{Y}(J\theta_L)(X),
\\[1ex]
g(UX,Y) =  \dh{}(J\theta_L)(X,Y).
\end{gather*}
We refer to \cite{VSC} for details about the vertical covariant
derivative operator $\DV{}$ and the horizontal exterior derivative
$\dh{}$ in these determining equations.

If we think of such a construction in the Finslerian case under
consideration, with $L=E$, it is natural to assume that $J$ is
homogeneous of degree zero and then $K$ has the same property. If
we then assume that this $K$ satisf\/ies the assumptions we needed
for the hierarchy of f\/irst integrals, namely $K=\ov{K}$ and the
condition \eqref{K}, then, unfortunately, there is still no
guarantee that $N_R=0$. Indeed, in the Finslerian case, although
the generically f\/ive dif\/ferent components of $N_R$ reduce to
two, and these considerably simplify further in view of the above
assumptions on $K$, in the end we run out of luck and they do not
automatically vanish! Note in passing that the vertical and
horizontal covariant derivative operators which appear in the
intrinsic characterization of the tensors $K$ and $U$ above are
associated, for the Finslerian case, to the Berwald connection. It
is very unlikely, however, that a dif\/ferent selection of
horizontal subspaces to describe these tensors (i.e.\ the choice
of one of the other linear connections, familiar in Finsler
geometry) might have any ef\/fect on the problem of involutivity
of the integrals. It further remains an open problem whether there
are dif\/ferent, more direct ways of obtaining some form of
involutivity. Also, the issue of functional independence of the
integrals, as treated in \cite{TopMat} for the Riemannian case,
has not been addressed.

\subsection*{Acknowledgement} This work has been partially
supported by the European Union through the FP6 Marie Curie RTN
{\em ENIGMA} (Contract number MRTN-CT-2004-5652).

\pdfbookmark[1]{References}{ref}
\LastPageEnding

\end{document}